\documentstyle[12pt]{article}
\newcommand{\setleftmargin}[1]{
	\addtolength{\textwidth}{\oddsidemargin}
	\addtolength{\textwidth}{1in}
	\addtolength{\textwidth}{-#1}
	\setlength{\oddsidemargin}{-1in}
	\addtolength{\oddsidemargin}{#1}
	\setlength{\evensidemargin}{\oddsidemargin}
}
\newcommand{\setrightmargin}[1]{
	\setlength{\textwidth}{8.5in}
	\addtolength{\textwidth}{-\oddsidemargin}
	\addtolength{\textwidth}{-1in}
	\addtolength{\textwidth}{-#1}
}
\newcommand{\settopmargin}[1]{
	\addtolength{\textheight}{\topmargin}
	\addtolength{\textheight}{1in}
	\addtolength{\textheight}{\headheight}
	\addtolength{\textheight}{\headsep}
	\addtolength{\textheight}{-#1}
	\setlength{\topmargin}{-1in}
	\addtolength{\topmargin}{-\headheight}
	\addtolength{\topmargin}{-\headsep}
	\addtolength{\topmargin}{#1}
}
\newcommand{\setbottommargin}[1]{
	\setlength{\textheight}{11in}
	\addtolength{\textheight}{-\topmargin}
	\addtolength{\textheight}{-1in}
	\addtolength{\textheight}{-\footskip}
	\addtolength{\textheight}{-#1}
}
\newcommand{\setallmargins}[1]{
	\settopmargin{#1}
	\setbottommargin{#1}
	\setleftmargin{#1}
	\setrightmargin{#1}
}
\setallmargins{1.2in}
\title{ 
The Generalized Witt Algebras using Additive Maps }
\author{ Ki-Bong Nam   and Mook Ok Wang}
\begin{document}
\maketitle
\begin{abstract}
Kawamoto generalized the Witt algebra using $F[x_1^{\pm 1},\cdots,x_n^{\pm 1}]$
instead of $F[x_1,\cdots,x_n].$
We construct the generalized Witt algebra $W(g_p,n)$ by using
an additive map $g_p$ from a set of integers into
a field of characteristic zero where $1\leq p \leq n.$
We show that the Lie algebra $W(g_p,n)$ is simple if
$g_p$ is injective, and also the Lie algebra
$W(g_p,n)$ has no ad-diagonalizable elements.
\end{abstract}

\newtheorem{lemma}{Lemma}
\newtheorem{prop}{Proposition}
\newtheorem{thm}{Theorem}
\newtheorem{coro}{Corollary}
\newtheorem{definition}{Definition}

\section{Introduction}

  Let $F$ be a field of characterisitic 
zero. 
The Witt algebra 
is called the general algebra by Rudakov [8].  
Kac [2] studied the generalized Witt algebra on 
the $F$-algebra in the formal power series 
$F[[x_1,\cdots ,x_n]]$ for a positive integer $n$.  
Nam
[5] constructs the Lie algebra on the $F$-subalgebra  
$F[e^{\pm x_1},\cdots ,e^{\pm x_n}, x_1,\cdots ,x_{n+m}]$ 
in the formal power series  
$F[[x_1,\cdots ,x_{n+}]]$ 
for the positive integers $n$ and $m.$

The Witt algebra $W(n)$ has a basis 
$$\{x_1^{a_1}\cdots x_n^{a_n}\partial_i |a_1,\cdots ,a_n\in N, 
1\leq i\leq n\}$$
with Lie bracket on basis elements
\begin{eqnarray*}
& &[x_1^{a_1}\cdots x_n^{a_n}\partial_i, 
x_1^{b_1}\cdots x_n^{b_n}\partial_j]\\ 
&=&b_ix_1^{a_1+b_1}\cdots x_n^{a_n+b_n}x_i^{-1}\partial_j 
-a_jx_1^{a_1+b_1}\cdots x_n^{a_n+b_n}x_j^{-1}\partial_i 
\end{eqnarray*}
where $N$ is the set of non-negative integers.

   Consider the generalized Witt algebra  $W(n,m)$ having a basis 
$$\{e^{a_1 x_1}\cdots e^{a_n x_n} x_1^{u_1}\cdots x_{n+m}^{u_{n+m}}\partial_k| 
a_1,\cdots ,a_n,u_1,\cdots ,u_{n+m}\in Z,1\leq k\leq n+m\}$$ 

with Lie bracket on basis elements given by 
$$[e^{a_1 x_1}\cdots e^{a_n x_n} x_1^{u_1}\cdots x_{n+m}^{u_{n+m}}\partial_i, 
e^{b_1 x_1}\cdots e^{b_n x_n} x_1^{t_1}\cdots x_{n+m}^{t_{n+m}}\partial_j]$$ 
$$=b_ie^{a_1x_1+b_1 x_1}
\cdots e^{a_nx_n+b_n x_n} x_1^{u_1+t_1}\cdots x_{n+m}^{u_{n+m}+t_{n+m}}
\partial_j$$ 
$$+t_ie^{a_1x_1+b_1 x_1}
\cdots e^{a_nx_n+b_n x_n} x_1^{u_1+t_1}\cdots x_{n+m}^{u_{n+m}+t_{n+m}}
x_i^{-1}\partial_j$$ 
$$-a_ie^{a_1x_1+b_1 x_1}
\cdots e^{a_nx_n+b_n x_n} x_1^{u_1+t_1}\cdots x_{n+m}^{u_{n+m}+t_{n+m}}
\partial_i$$ 
$$-u_ie^{a_1x_1+b_1 x_1}
\cdots e^{a_nx_n+b_n x_n} x_1^{u_1+t_1}\cdots x_{n+m}^{u_{n+m}+t_{n+m}}
x_j^{-1}\partial_i$$

where  $b_i=0$ if  $n+1\leq i \leq n+m,$
and $a_j=0$ if $n+1\leq j \leq n+m$ (see [4, 5, 7]).
In [5, 8], it is noted that 
the Lie subalgebra of $W(n,m)$ is the Witt algebra $W(n)$ on 
$F[x_1,\cdots ,x_m]$.
  Let $g_p$ be an additive map from $Z$ into $F$ where
$1\leq p \leq n$.  Let us define as $W(g_p,n)$ the Lie algebra with basis
\begin{eqnarray}\label{map1}
& &B:=\{{a_1 \choose i_1} \cdots {a_n \choose i_n}_k|a_1,\cdots ,a_n,i_1,\cdots ,i_n\in Z,
1\leq k \leq n\}            
\end{eqnarray}
and a Lie bracket on basis elements given by
\begin{eqnarray}\label{map10}
& &[{a_1 \choose i_1} \cdots {a_n \choose i_n}_k,
{b_1 \choose j_1} \cdots {b_n \choose j_n}_l]\\ \nonumber
&=&g_k(b_k){a_1+b_1 \choose i_1+j_1} \cdots {a_n+b_n \choose i_n+j_n}_l\\ \nonumber
&+&j_k{a_1+b_1 \choose i_1+j_1} \cdots {a_k+b_k \choose i_k+j_k-1} 
{a_{k+1} +b_{k+1} \choose i_{k+1}+j_{k+1}} \cdots {a_n+b_n \choose i_n+j_n}_l\\ \nonumber
&-&g_l(a_l){a_1+b_1 \choose i_1+j_1} \cdots {a_n+b_n \choose i_n+j_n}_l\\ \nonumber
&-&i_l{a_1+b_1 \choose i_1+j_1} \cdots 
{a_{l}+b_l \choose i_l+j_l-1} 
{a_{l+1}+b_{l+1} \choose i_{l+1}+j_{l+1}} 
\cdots {a_n+b_n \choose i_n+j_n}_k.
\end{eqnarray}
It follows from [1, 4, 5, 6] that 
the above bracket is extended linearly to the given basis $B$.  
Also, it is not hard to show 
that the above bracket satisfies the Jacobi identity.  In section 2,
we will prove the following main theorems.   Throughout 
the paper, a given map $g_p$ means an
additive and injective map from a set of integers 
into a field of characteristic zero.

\bigskip

\noindent
{\bf Theorem 1}  The Lie algebra $W(g_p,n)$ is simple.

\bigskip

\noindent
{\bf Theorem 2}  For any automorphism of $\theta  \in Aut (W(g,1))$ 
$$\theta ({0 \choose 1}_1 )=\sum_j C_j
({0 \choose j}_1 )$$
where $C_j\in F.$

\bigskip

\noindent
{\bf Theorem 3} Each derivation of $W(g,1)_+$ can be written
as a sum of an inner derivation and a scalar derivation \cite{Nam}.

\bigskip

  \section{ Simplicities of $W(g_p,n)$}

  The Lie algebra $W(g_p,n)$ has a $Z^n$-gradation as mentioned in 
[3] : that is, 
\begin{eqnarray}\label{map30}
& &W(g_p,n)=\bigoplus _{(a_1,\cdots ,a_n)\in Z^n} W_{(a_1,\cdots ,a_n)}
\end{eqnarray}
where each is the subspace of $W(g_p,n)$ with a basis 
$$\bar B=\{{a_1 \choose i_1}\cdots {a_n \choose i_n}_l|
a_1,\cdots,a_n,i_1,\cdots,i_n\in Z,1\leq l \leq n\}$$
 Let $W_{(a_1,\cdots ,a_n)}$ denote the $(a_1,\cdots ,a_n)$-homogeneous 
component of $W(g_p,n)$ and elements in   
the $W_{(a_1,\cdots ,a_n)}$-homogeneous $(a_1,\cdots ,a_n)$-homogeneous
elements.   Note that the 
${(0,\cdots ,0)}$
- homogeneous component 
is isomorphic to the Witt algebra $W(n)$ [8]. 
From now on let $(0,\cdots,0)$-homogeneous 
component denote the $0$-homogeneous component.
 For the simplicity of  $W(g_p,n)$, 
we define the map  as an additive and injective map. 
Now we introduce a lexicographic ordering of two basis 
elements of $W(g_p,n)$ as follows :
for any two elements 
${a_1 \choose i_1}\cdots {a_n \choose i_n}_l,$
${b_1 \choose j_1}\cdots {b_n \choose j_n}_k$
, we have  
$${a_1 \choose i_1}\cdots {a_n \choose i_n}_l>
{b_1 \choose j_1}\cdots {b_n \choose j_n}_k$$
 if $(a_1,\cdots ,a_n,i_1,\cdots ,i_n,l)
 >(b_1,\cdots ,b_n,j_1,\cdots ,j_n,k)$ by the natural
lexicographic ordering in $Z^{2n+1}$.

  For any element  $l\in W(g_p,n)$, 
$l$ can be written as follows
using the ordering and gradation
: 
$$l=\sum _{i_1,\cdots,i_n,p}C(i_1,\cdots ,i_n,p){a_{11} \choose i_1} 
\cdots   {a_{1n} \choose i_n}_l +\cdots $$
$$+\sum _{j_1,\cdots,j_n,q}C(j_1,\cdots ,j_n,q){a_{t1} \choose i_1} 
\cdots   {a_{tn} \choose j_n}_q$$

where $1\leq p,\cdots, q\leq n,$  $(a_{11},\cdots ,a_{1n}) 
>\cdots >(a_{11},\cdots ,a_{1n})$ 
and 
$$C(i_1,\cdots ,i_n,p),\cdots ,
C(j_1,\cdots ,j_n,q)\in F.$$

Next, define the string number $st(l)=t$ for $l$ (see [5, 6]), and $l_p(l)$ 
as $max \{i_1,\cdots ,i_n,\cdots ,j_1,\cdots ,j_n\}.$
  For any basis element  ${a_1 \choose i_1} \cdots {a_n \choose i_n}_l$ 
in $\bar B$, let us refer to $a_1,\cdots ,a_n$
 as upper indices  and $i_1,\cdots ,i_n$ as lower indices.

\bigskip

\noindent
{\bf Remark 2.1}  If $g_p$ is an inclusion, then $W(g_p,n)=W(n,0)$ 
is the generalized 
Witt algebra which is studied by Nam [5, 6].

  \begin{lemma}  If $l$ is any non-zero element, then the ideal $<l>$ 
generated by $l$ 
contains an element whose lower indices are positive.
\end{lemma}
{\it  Proof.}   Let $l$ be a nonzero element of $W(g_p,n)$.  
Take an element  $M={0 \choose j_1} \cdots {0 \choose j_n}_t$
such 
that $j_1>> \cdots >>j_n$ and $t$ such that either or $a_t\neq 0$ or
$i_t\neq 0$ in $l$, where $a>>b$ means $a$ is 
sufficiently larger than $b$.  Then $0\neq [M,l]$ 
is the required element.
\quad $\Box$
  \begin{lemma}  If an ideal $I$ of $W(g_p,n)$ contains ${0 \choose 0}\cdots {0 \choose 0}_l$ 
 for $1\leq l \leq n,$  then $I=W(g_p,n)$.
\end{lemma}
 {\it Proof.} Since $W(n)$ is a simple subalgebra of $W(g_p,n)$ 
with the basis 

\noindent
$\{ {0 \choose i_1}\cdots {0 \choose i_n}_l|i_1,\cdots ,i_n\in Z,
1\leq l \leq n\},$ the ideal
$< {0 \choose 0}\cdots {0 \choose 0}_i>$
contains $W(n)$, where 
$<{0 \choose 0}\cdots {0 \choose 0}>$ is the 
ideal of  which is generated by ${0 \choose 0}\cdots {0 \choose 0}_l$
for a fixed $t\in \{1,\cdots ,n\}$
 (see [8]).  It follows from Lemma 2 of [5] and Lemma 3 of [6] 
that the 
basis elements of $W(g_p,n)$ are contained in $I$
 by using the injectiveness of $g_p$. 
\quad $\Box$

\begin{thm}  The Lie algebra $W(g_p,n)$ is simple.
\end{thm}
 {\it Proof. } It is not difficult to prove 
this theorem by induction on $st(l)$ for any element $l$ in 
any ideal $I$ of $W(g_p,n)$ and by considering 
the one to one properties of $g_p$. 
\begin{coro}  The Lie algebra $W(n,0)$ is simple.
\end{coro}
 {\it Proof.}  If we take an additive 
embedding $g_p:Z\to F$, $1\leq p \leq n$, 
then we get the required result (see 
[5, 6]). 
\quad $\Box$

  It is an interesting problem to find all the automorphisms 
of the subalgebra 
$W(g,1)$ of $W(g_p,n)$.

\begin{thm}   For any automorphism $\theta \in Aut(g_p,n)$
$$\theta ({0 \choose 1}_1)=\sum _j C_j {0 \choose j }_1$$ 
where $C_j\in F$.
\end{thm}
  {\it Proof.}  It is not 
difficult to prove this theorem using the gradation of (2.1)   
and the $W(g,1)$
acting of  ${0 \choose 0}_1$ on the 
zero homogeneous component $W_0$ whose vector basis is $\{{0 \choose i}_1|i\in Z\}$
 as an ad-map.
\quad $\Box$

\bigskip

If we consider the Lie subalgebra of $W(g,1)$
such that all the lower indices are zero, then
this subalgebra is a Block algebra $W(1)$
which is also called the centerless Virasoro
algebra \cite{Blo}.
Thus all the automorphisms 
of this Lie algebra can be determined by Theorem 3 in [1].  

Consider the Lie subalgebra $W(1,1)^+$ of $W(n,m)$ in
\cite{Nam} with basis
$$\{e^{ax}x^iy^j\partial_x,
e^{bx}x^ly^m\partial_y| a,b\in Z, i,j,l,m\in N\}$$
where $N$ is the set of non-negative integers.

\bigskip
{\noindent}
{\bf Conjecture}
For any $\theta \in Aut (W(1,1)^+)$
$\theta (y\partial_y)=(\alpha y +\beta)\partial_y$
for some $0\neq \alpha, \beta\in F.$

\bigskip

The element $l\in W(g_p,n)$ is 
an ad-diagonalizable element if $[l,m]=\alpha (m) m$
 for any $m\in B$ given in (1) and for some $\alpha (m)\in F$.
\begin{prop}  The Lie algebra $W(g,n)$ 
has no ad-diagonalizable element with respect to 
the basis given in (\ref{map1}).
\end{prop}
 {\it Proof.} Since $W(g,n)$ is $Z^n$-graded Lie algebra 
all the ad-diagonalizable elements are in 
the $(0,\cdots,0)$-homogeneous component. 
$W_{(0,\cdots ,0)}$ is isomorphic to $W(n)$ as Lie algebras 
where $W(n)$ is the well known Witt algebra [8].  
Thus all the ad-semisimple elements of $W(g_p,n)$ 
are of the form  $\sum _{i=1}^n C_i {0 \choose 1 }$, 
where $C_i\in F$.  But   
 $[\sum _{i=1}^n C_i {0 \choose 1 }_i,
{a \choose 0}_j]\neq \alpha {a \choose 0}_j$
for any $\alpha \in F.$
  Therefore, we proved the proposition.
\quad $\Box$

\bigskip

\noindent
{\bf Remark 2.2}  For another proof of  
Proposition 2.6, see Corollary 1 of [5].

\bigskip

\section{Derivation of $W(g,1)_+$}

Consider the subalgebra $W(g,1)_+$ of $W(g_p,n)$ with basis
$$B_+:=\{{a\choose i}_1 |a\in Z,i\in N\}$$
where $N$ is the set of non-negative integers.

In this section we determine all the derivations of the
Lie algebra $W(g,1)_+$. Ikeda and Kawamoto found all
the derivations of the Kawamoto algebra $W(G,I)$ \cite{Kaw}
in their paper \cite{Ik}. It is very important to
find all the derivations of a given Lie algebra to
compare with other Lie algebras.
Let $L$ be a Lie algebra over any field $F.$
An $F$-linear map $D$ from $L$ to $L$ is a derivation
if $D([l_1,l_2])=[D(l_1),l_2]+[l_1,D(l_2)]$
for any $l_1,l_2\in L.$

Let $L$ be a Lie algebra over any field $F$. Define
the derivation $D$ of $L$ to be a scalar derivation if
for all basis elements $l$ of $L$,
$D(l)=f_l l$ for some scalar $f_l\in F$ \cite{Blo}, \cite{Nam1}.

We need the following lemma.
\begin{lemma}
Let $D$ be a derivation of $W(g,1)_+$. If $D(\partial)=0,$
then $D=f ad_{\partial} +S$ where
$f\in F$ and $S$ is a scalar derivation.
\end{lemma}
{\it Proof.}
It is not difficult to prove this lemma
using the gradation and the ordering of the Lie algebra
$W(g,1)_+.$ \cite{Nam}.
\quad $\Box$

\begin{thm}
Each derivation of $W(g,1)_+$ can be written
as a sum of an inner derivation and a scalar derivation.
\end{thm}
{\it Proof.}
Let $D$ be any derivation of $W(g,1)_+$. Then
$D(\partial)=f\partial$ for some
$f\in F[e^{\pm x},x]$. Since
$\partial : F[e^{\pm x},x] \to
F[e^{\pm x},x] $ is onto, there is a function
$g\in 
F[e^{\pm x},x] $ 
such that $\partial (g)=f.$ Then
$ad_{g\partial} (\partial)=[\partial,g\partial]=f\partial=D(\partial)$.
We have $(D-ad_{g\partial})(\partial)=0.$
By Lemma 3 we have $D=ad_{g\partial}+c ad_{\partial} +S$
\cite{Nam}. Therefore, we have proved the theorem.
\quad $\Box$

\bigskip

\noindent
{\bf Remark 3.1}
Let $Z^+$ and $F^+$ be the additive group of $Z$ and
$F$ respectively. For any $g\in Hom(Z^+,F^+),$ we have
$g(1)=m$ for some fixed $m\in F,$ thus the Witt algebra
$W(n)$ can not be changed by using the idea of the additive 
map in this paper \cite{Kaw1}.

\bigskip

{\it \small Department of Mathematics, UW-Madison, WI 53706}

{\it \small e-mail :nam@math.wisc.edu}

\bigskip

{\it \small Department of Mathematics, Hanyang University-Ansan, Ansan, Korea}
\end{document}